# Representation Theory Approach to the Polynomial Solutions of q - Difference Equations : $U_q(sl(3))$ and Beyond

**V.K. Dobrev**[*,*], **P. Truini**[†] and **L.C. Biedenharn**[♯]

[*] Arnold-Sommerfeld-Institut für Mathematische Physik,
Technische Universität Clausthal, Germany

[†] Istituto Nazionale di Fisica Nucleare, Sezione di Genova, Italy
Dipartimento di Fisica, dell'Universitá di Genova, Italy

[♯] Center for Particle Physics, Department of Physics,
University of Texas at Austin, USA

**Abstract**

A new approach to the theory of polynomial solutions of $q$–difference equations is proposed. The approach is based on the representation theory of simple Lie algebras $\mathcal{G}$ and their $q$ - deformations and is presented here for $U_q(sl(n))$. First a $q$–difference realization of $U_q(sl(n))$ in terms of $n(n-1)/2$ commuting variables and depending on $n-1$ complex representation parameters $r_i$, is constructed. From this realization lowest weight modules (LWM) are obtained which are studied in detail for the case $n=3$ (the well known $n=2$ case is also recovered). All reducible LWM are found and the polynomial bases of their invariant irreducible subrepresentations are explicitly given. This also gives a classification of the quasi-exactly solvable operators in the present setting. The invariant subspaces are obtained as solutions of certain invariant $q$-difference equations, i.e., these are kernels of invariant $q$-difference operators, which are also explicitly given. Such operators were not used until now in the theory of polynomial solutions. Finally the states in all subrepresentations are depicted graphically via the so called Newton diagrams.

[*] Permanent address: Institute for Nuclear Research and Nuclear Energy, Sofia, Bulgaria

## 1. Introduction.

There has been impressive progress in the last few years toward solving a fundamental problem posed by Bochner [1]: characterize differential operators possessing orthogonal polynomial eigenfunctions. This progress has taken many guises corresponding to several very different areas of research. In physics, the motivation has been to obtain solvable interaction Hamiltonians for bound state [2] and scattering problems [3] of the Schrödinger equation by using *Lie-algebraic techniques* (defined precisely below). Such techniques have re-appeared recently in the discovery by Turbiner, [4], Shifman [5], and others, of a new class of physically significant spectral problems, the so-called "quasi-exactly solvable problems", in which *part* of the (bound state) spectrum can be found algebraically. The basic structure of these related approaches has been recently clarified by the more mathematically oriented work of Kamran and Olver [6], [7], [8]. In classical applications of symmetry in quantum physics, using Lie group-theoretic methods, the Lie group appears as a symmetry group leaving the Hamiltonian invariant. A more general point-of-view arose in the late sixties with the introduction of spectrum generating groups [9] (and algebras) in which the symmetry group no longer necessarily left the Hamiltonian invariant but still retained some symmetry structure in the sense that the spectrum was "split but not mixed" in Gell-Mann's phrase [10]. In this more general approach, a second-order differential operator - - the Hamiltonian - - admits a spectrum generating symmetry if the Hamiltonian can be written as a quadratic in the generators (first-order differential operators) of a finite dimensional Lie algebra. The surprisingly successful interacting boson model of Iachello and Arima [2] is of this form with a compact Lie group symmetry ($SU(6)$); non-compact Lie group symmetries ($SL(3, {I\!\!R}), Sp(6, {I\!\!R}), \ldots$) have also enjoyed favor [11].

Following Kamran and Olver [7] we define a differential operator to be *Lie-algebraic* if it can be written as a quadratic combination of first-order differential operators (vector fields) which generate a finite-dimensional Lie algebra. A finite-dimensional Lie algebra is said to be [7] *quasi-exactly solvable* if it has a finite dimensional representation on some subspace of the space of smooth functions. In this terminology, a *quasi-exactly solvable differential operator* is a Lie-algebraic differential operator corresponding to a quasi-exactly solvable Lie algebra (the "hidden symmetry algebra"). As discussed by Shifman and Turbiner [5], the spectrum associated with the finite-dimensional representation spaces can be computed algebraically.

As one might expect, the one-dimensional case has been solved completely. It is useful to review these results here, since this case shows most clearly the essential underlying ideas: **Theorem** [7] (Gonzalez-Lopez, Kamran, Olver): *Every finite-dimensional quasi-exactly solvable Lie algebra of first order differential operators in one (real or complex)* *variable is locally equivalent to a subalgebra of one of the Lie algebras*

$$g_r \;=\; Span \;\left\{ \frac{d}{dz}, \; z\frac{d}{dz}, \; z^2\frac{d}{dz} - rz, \; 1 \right\} \;, \tag{1}$$

*where* $r \in \mathbb{Z}_+$ *is a non-negative integer. The associated representation module consists of polynomials of degree* $r$.

A cohomological interpretation of the integer $r$ has been given by Gonzalez-Lopez et al [8]. (The fact that $r$ must be integral has been called "quantized cohomology", and these authors [8] suggest the underlying reason as a problem for further investigation.)

The words "locally equivalent" in this theorem has a precise and important meaning: Two differential operators are said to be *locally equivalent* if they can be mapped into each other by a combination of change of independent variable: $x = f(z)$, and a (not necessarily unitary) "gauge transformation": $O \to e^{g(z)} O e^{-g(z)}$, where $O$ is a differential operator.

As an abstract Lie algebra $g_r$ in the result above is a central extension of the subalgebra:

$$h_r \;\equiv\; \left\{ J_0 = z\frac{d}{dz} - \frac{r}{2} \;,\; J_+ = rz - z^2\frac{d}{dz} \;,\; J_- = \frac{d}{dz} \;, \right\} \;. \tag{2}$$

It is sufficient, as Kamran and Olver remark, to concentrate attention - - without loss of generality - - on the algebra $h_r$ rather than on the full algebra $g_r$, since any Lie algebraic operator of $g_r$ is automatically a Lie-algebraic operator of $h_r$.

The Lie algebra given by (2) is very familiar to physicists: it is just the Lie algebra with the commutation relations:

$$[J_0, J_\pm] \;=\; \pm J_\pm \;,\quad [J_+, J_-] \;=\; 2J_0 \;, \tag{3}$$

which generate the quantal angular momentum group, $SU(2)$. The *realization* of this algebra (which can be expressed in terms of *one* pair $(a^+, a)$ of boson operators), is probably less generally familiar but equally important:

Remark: *The realization, (2), is the simplest example of a general algebraic technique for constructing <u>all</u> unitary irreps of all compact simple Lie groups, the algebraic Borel-Weil construction.*

The importance of this remark belies its ease of proof [12], [13], [14]. It is easily seen, for example, that the "quantization of cohomology" found in [8] is a direct consequence of the Borel-Weil construction (compactness) and as such appears to require no separate proof.

A further consequence of our remark is that the algebraic Borel-Weil construction implies a generalization of Lie algebraic quasi-exactly solvable Hamiltonians to encompass all simple Lie algebras. Whether or not this generalization encompasses *all* such structures, in the sense of the GKO theorem, is a question for further investigation.

The present paper is concerned with an application of our remark to the simple Lie algebra $A_\ell$, though it is clear, that the ideas are easily extendable for any simple Lie algebras $\mathcal{G}$. As a rank $\ell$ Lie algebra, the Cartan toroid must involve $\ell$ integers (rather than one as in the one-dimensional ($A_1$) algebra discussed above). Correspondingly, the algebraic Borel-Weil construction for $A_\ell$ must involve $\ell(\ell+1)/2$ complex variables (commuting bosons, which parametrize a Borel subgroup of $SL(\ell+1)$ [15]), thus going, even for $\ell = 2$, beyond the recent investigations of *two* variables [8].

We have chosen to discuss our extension of quasi-exactly solvable systems to more than two variables by using the more general *quantum group* $U_q(sl(n))$ instead of $sl(n)$, where $q$ is a generic parameter. This is technically a much more difficult construction, but it is at the same time much more general. Results for the Lie algebra $sl(n)$ are easily obtained by specializing our results to $q = 1$.

Note also that all our considerations are valid also for $U_q(sl(n, \mathbb{R}))$ if we restrict to real variables, real representation parameters and real $q$. Nothing is lost from the structure of the representations since all reducibilities occur when the representation parameters or their combinations are nonnegative integers.

We now give a summary of our approach and the results of the paper. The approach has several ingredients which are in fact available via the representation theory of simple Lie algebras $\mathcal{G}$ and their $q$ - deformations.

The first ingredient is the construction of a $q$ - difference realization of the representations of $U_q(\mathcal{G})$. For $U_q(sl(n))$ we use the construction of [15] for $U_q(su(n))$. The additional input w.r.t. [15] is to assign explicitly the eigenvalues of the Cartan generators. These are the representation parameters of our representations. Moreover, these eigenvalues are taken to be arbitrary complex numbers in order to obtain the most general representations. We also express all quantitites in terms of $n(n-1)/2$ commuting variables mentioned above ($n = \ell + 1$), and $q$ - difference and number operators in these variables.

Applying this realization to the constant function (independent of any of the variables) we obtain in fact realizations of lowest weight modules of $U_q(sl(n))$. In the present paper we consider the case when $q$ is not a nontrivial root of 1. Matters are arranged so that for generic representation parameters $r_i$, $i = 1, \ldots n-1$, these representations are irreducible and infinite-dimensional. If some representation parameters (or certain combinations thereof) are nonnegative integers the representation is reducible. The invariant irreducible subspace, or subrepresentation, is finite-dimensional iff all representation parameters are nonnegative integers: $r_i \in \mathbb{Z}_+$. In all cases the subrepresentations have polynomial bases.

Already at this point we can make a classification of all quasi-exactly solvable operators, which is equivalent to classifying the subrepresentations. This classification will depend entirely on the values of the representation parameters. Such a classification is also equivalent to the multiplet classification approach developed in [16]. In particular, for specific values of the representation parameters $r_i$, namely, only when some of them are zero, the invariant subspaces will depend on less variables; however, they can still be infinite-dimensional. As in the case of functions in one and two variables [17] all this can be visualized using the relevant Newton diagrams [18] which for the finite-dimensional representations are certain convex polyhedra in $\mathbb{R}_+^{n(n-1)/2}$. Moreover, it is natural from our point of view to introduce also infinite Newton diagrams corresponding to infinite-dimensional nontrivial subrepresentations.

We also discuss the $q$-difference operators which determine the $q$-difference equations whose solutions are the invariant subspaces. These are operators which arise when our representations are reducible and which intertwine these representations. We should stress that such operators, though widely used in many branches of mathematics and physics, are not used in the literature mentioned above. In general, such operators correspond to each positive root of the root system of $\mathcal{G}$ arising whenever the corresponding representation parameters are nonnegative integers. This situation is well understood in the $q = 1$ case [19], and the algebraic part of it is developed also for the $U_q(\mathcal{G})$ case [20].

The above approach is implemented fully in the present paper to the case $U_q(sl(3))$ (cf. Sections 3.3-3.6, 4.), while the general case is postponed to a sequel of the paper. (In passing, we also recover the well known $U_q(sl(2))$ case (cf. Section 3.2.)) For $U_q(sl(3))$ we give the basis of all invariant subspaces. The explicit polynomial basis (in three variables) needs hypergeometric functions for compact notation. Here there are three possible invariant operators of which at most two exist or are relevant for a particular invariant subspace. There are two generic situations, one of which contains, as a speaial case, the finite dimensional representations. The Newton diagrams are given in Section 4.

## 2. Procedure for the construction of the representations.

The procedure is iterative. In fact, we have to use also $U_q(gl(n))$. Let us introduce first some notation. The basic $q$-number notation is:

$$[a] = \frac{q^{a/2} - q^{-a/2}}{q^{1/2} - q^{-1/2}}, \tag{4}$$

and it will be used also for diagonal operators $H$ replacing $a$ in (4). Following [15] our representations will be given in terms of $n(n-1)/2$ variables. For our purposes we denote these variables by $z_i^k$, $2 \leq k \leq n$, $1 \leq i \leq k-1$. Next we introduce the number operator $N_i^k$ for the coordinate $z_i^k$, i.e., $N_i^k z_j^m = \delta_{mk}\delta_{ij}z_j^m$ and the $q$ - difference operators $D_i^k$, which admit a general definition on a larger domain than polynomials, but on polynomials are well defined as follows:

$$D_i^k = \frac{1}{z_i^k}[N_i^k]. \tag{5}$$

Further we note that the representations of $U_q(sl(n))$ will be characterized by $n-1$ complex parameters $r_k \in \mathbb{C}$, $1 \leq k \leq n-1$.

We rewrite formulae (5.3), (6.10) and (6.22) from [15] in the following way:

$$\Gamma_n(E_{ij}) = \Gamma_{n-1}(E_{ij}) \, q^{\frac{1}{4}(N_i^n - N_j^n)} + q^{\frac{1}{4}\Gamma_{n-1}(E_{jj} - E_{ii})} \, z_i^n D_j^n , \quad i < j < n, \tag{6a}$$

$$\Gamma_n(E_{ij}) = \Gamma_{n-1}(E_{ij}) \, q^{\frac{1}{4}(N_j^n - N_i^n)} + q^{\frac{1}{4}\Gamma_{n-1}(E_{ii} - E_{jj})} \, z_i^n D_j^n , \quad n > i > j, \tag{6b}$$

$$\Gamma_n(E_{ii}) = \Gamma_{n-1}(E_{ii}) + N_i^n , \quad i < n, \tag{6c}$$

$$\Gamma_n(E_{nn}) = \Gamma_{n-1}(E_{nn}) - \sum_{i=1}^{n-1} N_i^n , \tag{6d}$$

$$\Gamma_n(E_{ni}) = q^{\frac{1}{4}\left(\sum_{j=1}^{i-1}\Gamma_{n-1}(E_{jj}) - \sum_{j=i+1}^{n-1}\Gamma_{n-1}(E_{jj})\right)} D_i^n , \quad i < n, \tag{6e}$$

$$\Gamma_n(E_{in}) = q^{\alpha_{ii}^n} z_i^n \left[\Gamma_{n-1}(E_{nn}) - \Gamma_{n-1}(E_{ii}) - \sum_{k=1}^{n-1} N_k^n\right] -$$

$$- \sum_{\substack{j=1 \\ j \neq i}}^{n-1} q^{\alpha_{ij}^n} z_j^n \, \Gamma_{n-1}(E_{ij}), \quad i < n, \tag{6f}$$

where

$$q^{\alpha_{ij}^n} \equiv q^{-\frac{1}{4}\left(\sum_{k=1}^{j-1}\Gamma_{n-1}(E_{kk}) - \sum_{k=j+1}^{n-1}\Gamma_{n-1}(E_{kk})\right)} q^{\left(\frac{\pm}{0}\right)\left(\frac{1}{2}(\Gamma_{n-1}(E_{nn}) - \sum_{k=1}^{n-1} N_k^n) + \frac{3}{4}\right)} \times$$

$$\times q^{\left(\frac{\pm}{0}\right)\left(\frac{1}{4}(N_i^n + N_j^n)\right)} q^{\left(\frac{\pm}{0}\right)\frac{1}{2}\left(\sum_{\substack{i<k<j \text{ if } i<j \\ j<k<i \text{ if } j<i}} N_k^n\right)},$$

$$\left(\frac{\pm}{0}\right) = \begin{cases} + & \text{for } i < j, \\ - & \text{for } i > j, \\ 0 & \text{for } i = j. \end{cases} \tag{7}$$

$\Gamma_{n-1}(E_{ij})$ are defined at the previous step, except $\Gamma_{n-1}(E_{nn})$ which adds the representation parameter $r_{n-1}$ and is given by:

$$\Gamma_{n-1}(E_{nn}) = \sum_{k=0}^{n-1} r_k = r^{n-1} + r_0 \ . \tag{8}$$

The parameter $r_0$ represents the center of $U_q(gl(n))$ and is decoupled later. Note that $(6a-d)$, $(6e)$, $(6f)$, respectively, correspond to $(5.3a\text{-}d)$, $(6.10)$, $(6.22)$) of [15], respectively. The additional input with respect to [15] is : 1) notation - we put indices on $\Gamma$ corresponding to the case we consider - thus this is made an iterative procedure; 2) we give values to the Cartan generators which values are consistent with previous knowledge from representation theory; 3) we introduce $q$ - difference operators $D_i^n$ to replace $\bar{z}_i^n$; 4) we have done also something artificial by using $\Gamma_{n-1}(E_{nn})$ - this is given by the 'total number' $r^{n-1}$ (which for finite dimensional representations is the number of boxes of the Young tableaux minus $n-1$) plus the number $r_0$ representing the center of $U_q(gl(n))$.

Denoting $\mathcal{Z}^n = \sum_{i=1}^n E_{ii}$ we have:

$$\begin{aligned}\Gamma_n(\mathcal{Z}^n) &= \sum_{i=1}^n \Gamma_n(E_{ii}) = r^{n-1} + r_0 + \Gamma_{n-1}(\mathcal{Z}^{n-1}) = \\ &= \sum_{i=1}^{n-1} r^i + nr_0 = \sum_{i=0}^{n-1} (n-i)r_i \ .\end{aligned} \tag{9}$$

Thus, as expected, $\Gamma_n(\mathcal{Z}^n)$ is central. Then the generators $H_i^n = \Gamma_n(E_{ii} - E_{i+1,i+1})$, $1 \leq i < n$, $\Gamma_n(E_{ij})$, $i \neq j$, form a $q$ - difference operator realization of $U_q(sl(n))$.

It is straightforward to obtain the explicit expressions for $\Gamma_n(E_{ij})$. In particular, we have

$$\Gamma_n(E_{ii}) = \sum_{j=0}^{n-i-1} N_i^{n-j} - \sum_{j=1}^{i-1} N_j^i + \sum_{j=0}^{i-1} r_j \ , \tag{10}$$

with the usual convention that a sum is zero if the upper limit is smaller than the lower limit. From this we obtain the expressions for the Cartan generators $H_i^n$ (as defined above):

$$H_i^n = 2N_i^{i+1} - r_i + \sum_{j=0}^{n-i-2}(N_i^{n-j} - N_{i+1}^{n-j}) + \sum_{j=1}^{i-1}(N_j^{i+1} - N_j^i) \quad i < n \ . \tag{11}$$

Let us illustrate things for $n = 2, 3$.

For $n = 2$ we have (we use only $(6e, f)$, $(11)$) :

$$\Gamma_2(E_{12}) = z_1^2[r_1 - N_1^2] = x[r - N_x] \ , \tag{12a}$$
$$\Gamma_2(E_{21}) = D_1^2 = D_x \ , \tag{12b}$$
$$H_1^2 = 2N_x - r_1 \ , \tag{12c}$$

where we have denoted $z_1^2 = x$, $N_1^2 = N_x$. Note that we have obtained the well known [21] realization of the $U_q(sl(2))$ representations with $X^+ = \Gamma_2(E_{12})$, $X^- = \Gamma_2(E_{21})$, $H = H_1^2$, depending on the representation parameter $r_1$, ($r_0$ being cancelled as expected). For $q = 1$ this coincides with (2) setting $H = 2J_0$, $X^\pm = J_\pm$.

Next we take $n = 3$ setting $N_i^3 = N_i$, $D_i^3 = D_i$, $z_1^3 = z$, $z_2^3 = y$, $r = r^2 = r_1 + r_2$. Besides this we renormalize the generators $\Gamma_3(E_{13})$ and $\Gamma_3(E_{31})$ so that they obey the standard relations [22], [20] (these are different in [15], cf. (6.17), (6.20)):

$$\Gamma_3(E_{13}) = \Gamma_3(E_{12})\Gamma_3(E_{23}) - q^{1/2}\Gamma_3(E_{23})\Gamma_3(E_{12}), \tag{13a}$$

$$\Gamma_3(E_{31}) = \Gamma_3(E_{32})\Gamma_3(E_{21}) - q^{-1/2}\Gamma_3(E_{21})\Gamma_3(E_{32}). \tag{13b}$$

Thus, we have:

$$\Gamma_3(E_{12}) = \Gamma_2(E_{12})\, q^{\frac{1}{4}(N_1^3 - N_2^3)} + q^{\frac{1}{4}\Gamma_2(E_{22} - E_{11})}\, z_1^3 D_2^3 =$$
$$= x[r_1 - N_x]\, q^{\frac{1}{4}(N_z - N_y)} + zD_y q^{\frac{1}{4}(r_1 - 2N_x)}, \tag{14a}$$

$$\Gamma_3(E_{21}) = \Gamma_2(E_{21})\, q^{\frac{1}{4}(N_1^3 - N_2^3)} + q^{\frac{1}{4}\Gamma_2(E_{22} - E_{11})}\, z_2^3 D_1^3 =$$
$$= D_x\, q^{\frac{1}{4}(N_z - N_y)} + yD_z q^{\frac{1}{4}(r_1 - 2N_x)}, \tag{14b}$$

$$H_1^3 = 2N_1^2 - r_1 + N_1^3 - N_2^3 = 2N_x - r_1 + N_z - N_y, \tag{14c}$$

$$H_2^3 = 2N_2^3 - r_2 + N_1^3 - N_1^2 = 2N_y - r_2 + N_z - N_x, \tag{14d}$$

$$\Gamma_3(E_{31}) = q^{\frac{1}{4}(2\Gamma_3(E_{22}) - \Gamma_2(E_{22}))}\, D_1^3 = D_z q^{\frac{1}{4}(r_1 + r_0 - N_x + 2N_y)}, \tag{14e'}$$

$$\Gamma_3(E_{32}) = q^{\frac{1}{4}\Gamma_2(E_{11})}\, D_2^3 = D_y q^{\frac{1}{4}(r_0 + N_x)}, \tag{14e''}$$

$$\Gamma_3(E_{13}) = q^{\frac{1}{4}(\Gamma_2(E_{22}) - 2\Gamma_3(E_{22}))}\, z_1^3 \left[\Gamma_2(E_{33}) - \Gamma_2(E_{11}) - \sum_{k=1}^{2} N_k^3\right] -$$
$$- q^{-\frac{1}{2}\Gamma_3(E_{22})} q^{\alpha_{12}} z_2^3 \Gamma_2(E_{12}) =$$
$$= z[r - N_x - N_z - N_y] q^{\frac{1}{4}(N_x - r_1 - 2N_y - r_0)} -$$
$$- xy[r_1 - N_x] q^{\frac{1}{4}(2r_2 - r_0 + N_x - N_z - 3N_y + 1)} \tag{14f'}$$

$$\Gamma_3(E_{23}) = q^{-\frac{1}{4}\Gamma_2(E_{11})}\, z_2^3 \left[\Gamma_2(E_{33}) - \Gamma_2(E_{22}) - \sum_{k=1}^{2} N_k^3\right] - q^{\alpha_{21}} z_1^3 \Gamma_2(E_{21}) =$$
$$= y\, [r_2 + N_x - N_z - N_y]\, q^{-\frac{1}{4}(N_x + r_0)} -$$
$$- z\, D_x\, q^{-\frac{1}{4}(2r - r_1 + r_0 + N_x - N_z - N_y + 1)}. \tag{14f''}$$

We now rescale the generators $E_{3i}, E_{i3}$ ($i = 1, 2$) so as to absorb the parameter $r_0$. (Such a rescaling should be done also for the general $U_q(gl(n))$ case.) Thus the realization of $U_q(sl(3))$ depends only on the parameters $r_1, r_2$, as in the classical case. The latter is obtained from (14) by setting $q = 1$, ($\partial_a = \frac{\partial}{\partial a}$):

$$\Gamma_3(E_{12}) = x(r_1 - x\partial_x) + z\partial_y, \tag{15a}$$

$$\Gamma_3(E_{21}) = \partial_x + y\partial_z , \tag{15b}$$

$$H_1^3 = 2x\partial_x - y\partial_y + z\partial_z - r_1 , \tag{15c}$$

$$H_2^3 = -x\partial_x + 2y\partial_y + z\partial_z - r_2 , \tag{15d}$$

$$\Gamma_3(E_{31}) = \partial_z , \tag{15e'}$$

$$\Gamma_3(E_{32}) = \partial_y , \tag{15e''}$$

$$\Gamma_3(E_{13}) = z(r - x\partial_x - z\partial_z - y\partial_y) - yx(r_1 - x\partial_x) \tag{15f'}$$

$$\Gamma_3(E_{23}) = y(r_2 + x\partial_x - z\partial_z - y\partial_y) - z\partial_x . \tag{15f''}$$

Of course, (15) can be derived from the induced representation picture [19].

### 3. Reducibility of the representations and invariant subspaces.

**3.1.** Let us apply the realization (6) to the function 1. Using the fact that $N_i^n \, 1 = 0 = D_i^n \, 1$ we have:

$$\Gamma_n(E_{ii}) \, 1 = \sum_{j=0}^{i-1} r_j = r^{i-1} , \quad H_i^n \, 1 = -r_i , \quad i \le n , \tag{16a}$$

$$\Gamma_n(E_{ni}) \, 1 = 0 , \quad i < n , \tag{16b}$$

$$\Gamma_n(E_{ij}) \, 1 = \Gamma_{n-1}(E_{ij}) \, 1 = \cdots = \Gamma_i(E_{ij}) \, 1 = 0 , \quad j < i < n , \tag{16c}$$

$$\Gamma_n(E_{ij}) \, 1 = \Gamma_{n-1}(E_{ij}) \, 1 = \cdots = \Gamma_j(E_{ij}) \, 1 , \quad i < j < n , \tag{16d}$$

$$\Gamma_n(E_{in}) \, 1 = q^{\frac{1}{4}\left(\sum_{k=i+1}^{n-1} r^{k-1} - \sum_{k=1}^{i-1} r^{k-1}\right)} z_i^n \left[r_i + \cdots + r_{n-1}\right] -$$
$$- \sum_{s=i+1}^{n-1} q^{\alpha_{is}^n} z_s^n \, \Gamma_s(E_{is}) \, 1 , \quad i < n , \tag{16e}$$

where

$$q^{\alpha_{is}^n} = q^{\frac{1}{4}\left(\sum_{k=s+1}^{n-1} r^{k-1} - \sum_{k=1}^{s-1} r^{k-1}\right)} q^{\frac{1}{2}(r^{n-1} - \sum_{k=1}^{n-1} N_k^n) + \frac{3}{4}} \times$$
$$\times q^{\frac{1}{4}(N_i^n + N_s^n)} q^{\frac{1}{2}\sum_{k=i+1}^{k=s-1} N_k^n} , \quad i < s . \tag{17}$$

It is straightforward to obtain the explicit expressions for $\Gamma_n(E_{ij}) \, 1$, $i < j \le n$, applying recursively (16d, e). In particular, we have:

$$\Gamma_n(E_{i,i+1}) \, 1 = \Gamma_{i+1}(E_{i,i+1}) \, 1 = q^{-\frac{1}{4}\sum_{k=1}^{i-1} r^{k-1}} [r_i] \, z_i^{i+1} . \tag{18}$$

Thus, we have obtained a lowest weight module (LWM) with lowest weight vector 1, (it is annihilated by the lowering generators $\Gamma_n(E_{ij})$, $j < i \le n$), and lowest weight $\Lambda$ such that $\Lambda(H_i) = -r_i$, (cf. (16a)). Generically this LWM is irreducible and then it is isomorphic to the Verma module with this lowest weight. The states in it correspond to the monomials of the Poincaré-Birkhoff-Witt basis of $U_q(\mathcal{G}^+)$, where $\mathcal{G}^+$ is the subalgebra of the raising generators. This is isomorphic to the monomials in the variables $z_i^k$. When the representation parameters $r_i$ or certain combinations thereof are nonnegative integers our representations are reducible. The full analysis of this is postponed to the sequel of this paper. Below we consider in detail the cases $n = 2$ and $n = 3$.

**3.2.** We start with $n = 2$ (though this example is well known). Let us apply $H = \Gamma_2(E_{11} - E_{22})$, $X^+ = \Gamma_2(E_{12})$, $X^- = \Gamma_2(E_{21})$ to the function 1. We use the fact that $N_x 1 = 0 = D_x 1$. Thus:

$$H\,1 \;=\; -r\,,\qquad X^+\,1 \;=\; x[r]\,,\qquad X^-\,1 \;=\; 0\,. \tag{19}$$

Thus, we obtain a lowest weight module with lowest weight vector 1 and lowest weight $\Lambda$ such that $\Lambda(H) = -r$. All states are given by powers of $x$, i.e., the basis is $x^k$ with $k \in \mathbb{Z}_+$ and the representation is infinite dimensional. The action of $U_q(sl(2))$ is given by:

$$X^+ x^k \;=\; [r-k] x^{k+1}\,,\qquad X^- x^k \;=\; [k] x^{k-1}\,,\qquad H x^k \;=\; (2k-r) x^k\,. \tag{20}$$

Clearly, if $r \notin \mathbb{Z}_+$ this representation is irreducible. Furthermore all states may be obtained by the application of $X^+$ to the LWV, i.e., :

$$(X^+)^k\,1 \;=\; x^k [r][r-1]\ldots[r-k+1]\,,\qquad k \in \mathbb{Z}_+\,. \tag{21}$$

Let $r \in \mathbb{Z}_+$, then $(X^+)^{r+1}\,1 = X^+ x^r [r]! = 0$. Thus, the states $x^k$ with $k = 0, 1, \ldots r$ form a finite-dimensional subrepresentation with $\dim = r + 1$. Note that the complement of this subrepresentation, i.e., the states $x^k$ with $k > r$, is not an invariant subspace.

Clearly, any polynomial in $H, X^\pm$, will preserve this invariant subspace and thus would be a quasi-exactly solvable operator.

The invariant subspace may be obtained as the solution of either one the following equations:

$$(X^+)^{r+1}\,f(x) \;=\; 0\,, \tag{22a}$$

$$(X^-)^{r+1}\,f(x) \;=\; 0\,, \tag{22b}$$

in the space of formal power series $f(x) = \sum_{k \in \mathbb{Z}_+} \mu_k x^k$. Note, however, that only (22b) (which is enough) was expected - this is an artefact of $n = 2$ simplifications. Indeed, only the operator in (22b) has the intertwining property (as in the classical case [19]):

$$(X^-)^{r+1}\,\Gamma_2(X)_r \;=\; \Gamma_2(X)_{r'}\,(X^-)^{r+1}\,,\qquad r' = -r - 2\,, \tag{23}$$

where $X = H, X^{\pm}$, and $\Gamma_2(X)_r$ is from (12) with explicit notation for the representation parameter of the two representations which are intertwined.

**3.3.** To the end of Section 3. we discuss $U_q(sl(3))$.

**3.3.1.** Let us apply (14) to the function 1 :

$$\Gamma_3(E_{12})\, 1 = x[r_1]\,, \tag{24a}$$
$$\Gamma_3(E_{21})\, 1 = 0\,, \tag{24b}$$
$$H_1\, 1 = -r_1\,, \tag{24c}$$
$$H_2\, 1 = -r_2\,, \tag{24d}$$
$$\Gamma_3(E_{31})\, 1 = 0\,, \tag{24e'}$$
$$\Gamma_3(E_{32})\, 1 = 0\,, \tag{24e''}$$
$$\Gamma_3(E_{13})\, 1 = q^{-\frac{1}{4}r_1} z[r] - q^{\frac{1}{4}(2r_2+1)} yx[r_1] \tag{24f'}$$
$$\Gamma_3(E_{23})\, 1 = y[r_2]\,, \tag{24f''}$$

Thus, we obtain a lowest weight module with lowest weight vector 1 and lowest weight $\Lambda$ such that $\Lambda(H_k) = -r_k$. All states are given by powers of $x, y, z$, i.e., the basis is generated by $x^j z^k y^\ell$ with $j, k, \ell \in \mathbb{Z}_+$. The action of $U_q(sl(3))$ is given by:

$$\Gamma_3(E_{12})\, x^j z^k y^\ell = [r_1 - j]\, q^{\frac{1}{4}(k-\ell)}\, x^{j+1} z^k y^\ell + [\ell] q^{\frac{1}{4}(r_1 - 2j)}\, x^j z^{k+1} y^{\ell-1}\,, \tag{25a}$$
$$\Gamma_3(E_{21})\, x^j z^k y^\ell = [j]\, q^{\frac{1}{4}(k-\ell)}\, x^{j-1} z^k y^\ell + [k] q^{\frac{1}{4}(r_1 - 2j)}\, x^j z^{k-1} y^{\ell+1}\,, \tag{25b}$$
$$H_1\, x^j z^k y^\ell = (-r_1 + 2j - \ell + k)\, x^j z^k y^\ell\,, \tag{25c}$$
$$H_2\, x^j z^k y^\ell = (-r_2 - j + 2\ell + k)\, x^j z^k y^\ell\,, \tag{25d}$$
$$\Gamma_3(E_{31})\, x^j z^k y^\ell = [k] q^{\frac{1}{4}(r_1 - j + 2\ell)}\, x^j z^{k-1} y^\ell\,, \tag{25e'}$$
$$\Gamma_3(E_{32})\, x^j z^k y^\ell = [\ell] q^{\frac{1}{4} j}\, x^j z^k y^{\ell-1}\,, \tag{25e''}$$
$$\Gamma_3(E_{13})\, x^j z^k y^\ell = q^{\frac{1}{4}(j - r_1 - 2\ell)}[r - j - k - \ell]\, x^j z^{k+1} y^\ell - $$
$$- q^{\frac{1}{4}(2r_2 + j - k - 3\ell + 1)}[r_1 - j]\, x^{j+1} z^k y^{\ell+1}\,, \tag{25f'}$$
$$\Gamma_3(E_{23})\, x^j z^k y^\ell = q^{-\frac{1}{4} j}\, [r_2 + j - k - \ell]\, x^j z^k y^{\ell+1} - $$
$$- [j] q^{-\frac{1}{4}(r_1 + 2r_2 + j - k - \ell + 1)}\, x^{j-1} z^{k+1} y^\ell\,. \tag{25f''}$$

In this Section we show the following results which parallel the classical situation (cf. [19]):

**1.** If $r_1$, or $r_2$, or $r+1 \in \mathbb{Z}_+$ this representation is reducible. It contains an irreducible subrepresentation which is infinite-dimensional, except when both $r_1$, $r_2 \in \mathbb{Z}_+$ ;

**2.** If $r_1$, $r_2$, $r + 1 \notin \mathbb{Z}_+$ this representation is irreducible and infinite dimensional.

**3.3.2.**   Clearly, if $r_1 \in \mathbb{Z}_+$ the representation (25) becomes reducible : the monomials $x^j z^k y^\ell$ with $j \leq r_1$ form an invariant subspace since from $(25a, f')$ we have:

$$\Gamma_3(E_{12}) \, x^{r_1} z^k y^\ell = [\ell] q^{-\frac{1}{4}r_1} \, x^{r_1} z^{k+1} y^{\ell-1} \,, \tag{26}$$

$$\Gamma_3(E_{13}) \, x^{r_1} z^k y^\ell = [r_2 - k - \ell] q^{-\frac{1}{2}\ell} \, x^{r_1} z^{k+1} y^\ell \,, \tag{27}$$

and all other operators are either lowering or preserving the powers of $x$. This invariant subspace may be described as the solution of the following $q$ - difference equation:

$$(D_x)^{r_1+1} \, f(x,y,z) = 0 \,. \tag{28}$$

Note that the operator in (28) has the intertwining property (as in the classical case [19]):

$$(D_x)^{r_1+1} \, \Gamma_3(X)_{r_1,r_2} = \Gamma_3(X)_{r'_1,r'_2} \, (D_x)^{r_1+1} \,, \qquad r'_1 = -r_1 - 2 \,, \ r'_2 = r+1 \,, \tag{29}$$

where $X = E_{ii} - E_{i+1,i+1}, E_{ij}, i \neq j$, $\Gamma_3(X)_{r_1,r_2}$ is taken from (14) with explicit dependence of the representation parameters of the two representations which are intertwined.

The subrepresentation obtained is infinite dimensional if $r_2 \notin \mathbb{Z}_+$ since the powers of $y, z$ are still unrestricted by $(25 f', f'')$.

**3.3.3.**   If $r_2 \in \mathbb{Z}_+$ the representation in (25) becomes reducible. In the classical case $(q = 1)$ the equation which singles out the invariant subspace is [19]:

$$\left( x \partial_z + \partial_y \right)^{r_2+1} \, f(x,y,z) = 0 \,, \quad q = 1 \,. \tag{30}$$

For the quantum case we have the following expression :

$$_q\mathcal{D}_2^{r_2+1} \, f(x,y,z) = 0 \,, \tag{31a}$$

$$_q\mathcal{D}_2^{r_2+1} = \sum_{s=0}^{k} \binom{k}{s}_q x^{k-s} \, D_z^{k-s} \, D_y^s \, q^{\frac{1}{4}s(N_z - r_1) + \frac{1}{4}(s-k)N_y - \frac{1}{4}kN_x} \,, \tag{31b}$$

which coincides with (30) for $q = 1$. The invariant subspace is infinite dimensional if $r_1 \notin \mathbb{Z}_+$.

As in the classical case [19] the explicit form (31b) of this operator may be checked by the intertwining property:

$$_q\mathcal{D}_2^{r_2+1} \, \Gamma_3(X)_{r_1,r_2} = \Gamma_3(X)_{r'_1,r'_2} \, _q\mathcal{D}_2^{r_2+1} \,, \qquad r'_1 = r+1 \,, \ r'_2 = -r_2 - 2 \,, \tag{32}$$

where $X = E_{ij}, i \neq j, E_{ii} - E_{i+1,i+1}$,   $\Gamma_3(X)_{r_1,r_2}$ is from (14) as in (29).

**3.4.** In this subsection we consider the case when both $r_k \in \mathbb{Z}_+$ .

**3.4.1.** Let $r_k \in \mathbb{Z}_+$ , $k = 1, 2$. Then there is a finite dimensional irreducible subspace of dimension:
$$d_{r_1,r_2} = \frac{1}{2}(r_1 + 1)(r_2 + 1)(r + 2) . \tag{33}$$

Thus, we recover the complete list of the finite dimensional irreps of $U_q(sl(3))$ and $SL(3)$, and by default, also the complete list of the finite dimensional unitary irreps of $U_q(su(3))$ and $SU(3)$ (we have assumed that $q$ is not a nontrivial root of 1). Below we give explicitly the basis of this subspace.

**3.4.2.** Let us consider first the special case $r_1 = 0$, $r_2 = r \in \mathbb{Z}_+$. (This is the case given in [17] for $q = 1$.) Then there is no $x$ dependence in the irreducible subrepresentation as we argued above. Since $\Gamma_3(E_{12}) \, 1 = 0$, all states can be built in the following way. Apply first $\Gamma_3(E_{23})$ to the function 1:

$$(\Gamma_3(E_{23}))^\ell \, 1 = [r][r-1]\ldots[r-\ell+1] \, y^\ell \tag{34}$$

thus we obtain $r + 1$ states $((\Gamma_3(E_{23}))^{r+1} \, 1 = 0)$. Then we apply the operator $\Gamma_3(E_{13})$ to each of these states:

$$(\Gamma_3(E_{13}))^k y^\ell = q^{-\frac{1}{2}k\ell} [r-\ell][r-\ell-1]\ldots[r-\ell-k+1] \, y^\ell z^k \tag{35}$$

thus for each fixed $\ell$ we obtain $r + 1 - \ell$ states.

Thus the number of states is (cf. (33)) :

$$\sum_{\ell=0}^{r} \sum_{k=0}^{r-\ell} 1 = \frac{1}{2}(r+1)(r+2) = d_{0,r} . \tag{36}$$

**3.4.3.** In a similar way we obtain the case $r_1 > 0$, however, hypergeometric finctions enter the game. Apply first $\Gamma_3(E_{12})$ to the function 1:

$$(\Gamma_3(E_{12}))^j \, 1 = [r_1][r_1-1]\ldots[r_1-j+1] \, x^j \tag{37}$$

thus we obtain $r_1 + 1$ states $((\Gamma_3(E_{12}))^{r_1+1} \, 1 = 0)$. Then we apply the operator $\Gamma_3(E_{13})$ to each of these states:

$$(\Gamma_3(E_{13}))^k x^j = \sum_{s=0}^{k} (-1)^s \beta_{ks} \, x^{j+s} z^{k-s} y^s . \tag{38}$$

Applying $\Gamma_3(E_{13})$ once more we obtain a recursion relation:

$$\beta_{k+1,s} = q^{\frac{1}{4}(j-r_1-s)} [r-j-s-k]\beta_{ks} + q^{\frac{1}{4}(2r_2+j-k-s+2)} [r_1-j-s+1]\beta_{k,s-1} \tag{39}$$

with the convention that $\beta_{k,k+1} = \beta_{k,-1} = 0$. Solving this, we have :

$$\beta_{k,s} = q^{\frac{1}{4}k(j-r_1) + \frac{1}{4}s(2r-r_1+2-k)} \binom{k}{s}_q \frac{[r_1-j]![r-j-s]!}{[r_1-j-s]![r-j-k]!} . \tag{40}$$

Note that:
$$\beta_{k,s} = 0 , \quad \text{if} \quad \begin{cases} s > \min(k, r_1-j) & \text{for } r \neq r_1 \\ s > k & \text{for } r = r_1 \end{cases} \tag{41}$$

Thus we have:

$$(\Gamma_3(E_{13}))^k \, \Gamma_3(E_{12}))^j \, 1 \; = \; x^j z^k q^{\frac{1}{4}k(j-r_1)} \frac{[r_1]!}{[r-j-k]!} \sum_{s=0} \binom{k}{s}_q \times$$
$$\times \frac{[r-j-s]!}{[r_1-j-s]!} q^{\frac{1}{4}s(2r-r_1+2-k)} \left(-\frac{xy}{z}\right)^s =$$

$$\tag{42}$$

$$= \begin{cases} x^j z^k q^{\frac{1}{4}k(j-r_1)} \frac{[r-j]![r_1]!}{[r-j-k]![r_1-j]!} \times & \\ \qquad \times \, {}_2F_1^q(-k, j-r_1; j-r; q^{\frac{1}{4}(2r-r_1+2-k)} \frac{xy}{z}) & \text{for } r \neq r_1 \\ x^j z^k q^{\frac{1}{4}k(j-r)} \frac{[r]!}{[r-j-k]!} \left(1 - q^{\frac{1}{4}(r+2-k)} \frac{xy}{z}\right)_q^k , & \text{for } r = r_1 \end{cases}$$

where ${}_2F_1^q$ is the $q$-hypergeometric function, $(1+w)_q^k$ is the $q$-binomial. From (42) we have:

$$(\Gamma_3(E_{13}))^k x^j \; = \; 0 , \quad k > r-j . \tag{43}$$

Thus, we have at this stage $(r_1+1)(2r+2-r_1)/2$ states. Now it remains to apply $(\Gamma_3(E_{23}))^\ell$ to $(\Gamma_3(E_{13}))^k \, \Gamma_3(E_{12}))^j \, 1$ to obtain the remaining $r_2(r_1+1)(r+1)/2$ states.

Let us denote $v_{kj} = \Gamma_3(E_{13}))^k \, \Gamma_3(E_{12}))^j \, 1$. If $r = r_1$, $(r_2 = 0)$, we have:

$$\Gamma_3(E_{23}) \, v_{kj} \; = \; -[j] x^{j-1} z^{k+1} q^{\frac{1}{4}((k-1)(r+1) - j(k+1))} \frac{[r]!}{[r-j-k]!} \times$$
$$\times \left(1 - q^{\frac{1}{4}(r+1-k)} \frac{xy}{z}\right)_q^{k+1} = \tag{44}$$
$$= -q^{-\frac{1}{2}(r+1)} [j] \, v_{k+1,j-1} ,$$

i.e., $\Gamma_3(E_{23})$ just transforms the states obtained so far. Thus, the number of the states is :

$$\sum_{j=0}^{r} \sum_{k=0}^{r-j} 1 \; = \; \frac{1}{2}(r+1)(r+2) \; = \; d_{r,0} \tag{45}$$

as expected (cf. (33)).

**Statement.** This subspace may be characterized as the solution space of two equations, namely, (28) and:

$$\left( xD_z \; q^{-\frac{1}{4}(N_y+N_x)} + D_y q^{\frac{1}{4}(N_z-N_x-r)} \right) f(x,y,z) \;=\; 0 , \tag{46}$$

which is (31) for $r_2 = 0$.

**Proof.** Let us write $f(x,y,z) = \sum_{j,k,\ell \in \mathbb{Z}_+} \lambda_{jk\ell} x^j z^k y^\ell$ and apply (46). Below we consider $\lambda_{jk\ell}$ for $j, k, \ell \in \mathbb{Z}$ with the convention that $\lambda_{jk\ell} = 0$ if $j < 0$, or $k < 0$, or $\ell < 0$. We obtain the following recursion relation for $\lambda_{jk\ell}$:

$$[k+1]\lambda_{j-1,k+1,\ell} + [\ell+1]\lambda_{jk,\ell+1} q^{\frac{1}{4}(k+\ell-r)} \;=\; 0. \tag{47}$$

We note that from (28) follows that $\lambda_{jk\ell} = 0$ if $j > r_1 = r$. From this and (47) follows that $\lambda_{jk\ell} = 0$ if $j + k > r$, (since $\lambda_{jk\ell} \sim \lambda_{r+1,k+j-r-1,\ell-j+r+1} = 0$ with proportionality coefficient being non-zero and finite for the hypothesis in consideration). Next we note that $\lambda_{jk\ell} = 0$ if $\ell > j$, (since $\lambda_{jk\ell} \sim \lambda_{-1,k+j+1,\ell-j-1} = 0$). We also note that $\lambda_{jk\ell}$ are proportinal to each other for $k + \ell$ fixed. In order to ensure the above properties automatically we are prompted to make the substitution:

$$j \to j + s, \quad k \to k - s, \quad \ell \to s, \quad \text{with } s \in \mathbb{Z}_+, \quad s \leq k. \tag{48}$$

Solving (47) with this we get:

$$\lambda_{j+s,k-s,s} \;=\; (-1)^s \binom{k}{s}_q q^{\frac{s}{4}(r+2-k)} \lambda_{jk0} , \tag{49}$$

where $\lambda_{jk0}$ are only subject to the vanishing condition for $j + k > r$. The latter is achieved if we suppose that $\lambda_{jk0} \sim 1/[r-j-k]_q!$. Thus we see that up to normalization $\lambda_{j+s,k-s,s}$ coincide with $\beta_{ks}$ for $r = r_1$. •

**3.4.4.** Thus, it remains to treat the case when $r_k \in \mathbb{N}$. First we note another case where the $\Gamma_3(E_{23})$ transforms the states:

$$\Gamma_3(E_{23}) \; v_{r-j,j} \;=\; -q^{-\frac{1}{2}(r+1)} [j] \; v_{r-j+1,j-1} , \tag{50}$$

Next we shall use the following general formula valid for arbitrary $r_k$:

$$v_{\ell k j} \equiv \Gamma_3(E_{23})^\ell \, \Gamma_3(E_{13})^k \, \Gamma_3(E_{12})^j \, 1 \; =$$

$$= \sum_{s=0}^{k} \sum_{n=0}^{\ell} (-1)^{s-n} \binom{k}{s}_q \binom{\ell}{n}_q q^{\frac{1}{4}\{(j-r_1)k - \ell j + (s-n)(r_1 + 2r_2 - k - \ell + 2)\}} \times$$

$$\times \frac{\Gamma_q(r_1+1)\Gamma_q(r-j-s+1)}{\Gamma_q(r_1-j-s+1)\Gamma_q(r-j-k+1)} \frac{\Gamma_q(r_2+j+s-k-n+1)}{\Gamma_q(r_2+j+s-k-\ell+1)} \frac{[j+s]!}{[j+s-n]!} \times$$

$$\times x^{j+s-n} \, z^{k-s+n} \, y^{\ell+s-n} \; . \tag{51}$$

Note that for $r_2 \in \mathbb{Z}$ the ratio $\Gamma_q(r_2+j+s-k-n+1)/\Gamma_q(r_2+j+s-k-\ell+1)$ means just $[r_2+j+s-k-n][r_2+j+s-k-n-1]\ldots[r_2+j+s-k-l+1]$.

A basis for our representation space is given by $v_{\ell k j}$ iff $\ell + k + j \leq r$, $0 \leq j \leq r_1$, $0 \leq \ell \leq r_2$. To show this it is enough to demonstrate that $v_{\ell k j}$ with any of the conditions violated is a linear combination of the representation states (which combination may be zero) - cf. e.g., (44) and (50). In particular, by induction, it is enough to show this for the boundary cases: $\ell = r_2 + 1$, $\ell + j + k = r + 1$. We give explicitly the boundary case $\ell = r_2 + 1$. Let us introduce the notation

$$A = \Gamma_3(E_{23}) \,, \quad B = \Gamma_3(E_{12}) \,, \quad C = \Gamma_3(E_{13}) \,. \tag{52}$$

Using the commutation relations it is easy to show that:

$$A^\ell C^k = q^{\frac{k\ell}{2}} C^k A^\ell \,, \qquad C^k B^j = q^{\frac{jk}{2}} B^j C^k \,, \tag{53}$$

$$B^j A^\ell = \sum_{t=0}^{\lambda} \binom{\ell}{t}_q \binom{j}{t}_q [t]! \, q^{\frac{(\ell-t)(j-t)}{2}} A^{\ell-t} C^t B^{j-t} \,, \quad \lambda = \min\{\ell, j\} \,. \tag{54}$$

The proof of the result stated above now follows very easily from (53), (54) and the fact that $A^{r_2+1} \, 1 = 0$. In fact

$$A^{r_2+1} C^k B^j \, 1 \; = \; q^{\frac{(r_2+1)k}{2}} \, C^k \times$$

$$\times \left( B^j A^{r_2+1} - \sum_{t=1}^{\lambda} \binom{r_2+1}{t}_q \binom{j}{t}_q [t]! \, q^{\frac{(r_2+1-t)(j-t)}{2}} A^{r_2+1-t} C^t B^{j-t} \right) 1 \,, \tag{55}$$

here $\lambda = \min\{r_2+1, j\}$. The first term in the parenthesis gives zero when applied to 1, while the second term yields a linear combination of states of the required form ($\ell \leq r_2$). This ends the proof. The boundary case $\ell + j + k = r + 1$ may be treated in a similar way.

**3.5.** Next, we discuss the case when $r+1 \in \mathbb{Z}_+$, but $r_k \notin \mathbb{Z}_+$. In the classical case there is a nontrivial invariant subspace singled out by the following equation [19]:

$$\mathcal{D}_3^{r+2} f(x,y,z) = 0, \tag{56a}$$

$$\mathcal{D}_3^{r+2} = \sum_{s=0}^{r+2} \frac{(-1)^s}{(r_1+1-s)} \binom{r+2}{s} (\hat{X}_1^+)^{r+2-s} (\hat{X}_2^+)^{r+2} (\hat{X}_1^+)^s, \tag{56b}$$

$$\hat{X}_1^+ = \partial_x, \qquad \hat{X}_2^+ = x\partial_z + \partial_y. \tag{56c}$$

For this equation one uses the singular vector in (8.41) of [19]. The necessary singular vector in the $q$ case is obtained from the classical case by just replacing the numbers with $q$-numbers [20]. However, in order to generalize this to the $q$ - case we first rewrite (56b) with the explicit substitution of (56c). We obtain:

$$\mathcal{D}_3^{r+2} = -\sum_{s=0}^{r+2} \frac{(r+2)!\,\Gamma(-1-r_1)}{\Gamma(-r_1+r+2-s)} \binom{r+2}{s} \partial_x^s (x\partial_z + \partial_y)^{r+2-s} \partial_x^{r+2-s} = \tag{57a}$$

$$= -\sum_{s=0}\sum_{t=0} \frac{(r+2)!^2\,\Gamma(-1-r_1)}{s!\,t!\,(r+2-s-t)!\,\Gamma(r-r_1+2-s)} \times$$

$$\times \partial_z^{r+2-t}\,\partial_y^t\,\partial_x^t \prod_{u=1}^{r+2-s-t} (x\partial_x - t + 1 - u). \tag{57b}$$

Now we have the following expression for the $q$ - difference equation in our case:

$$_q\mathcal{D}_3^{r+2} f(x,y,z) = 0, \tag{58a}$$

$$_q\mathcal{D}_3^{r+2} = \sum_{s=0}\sum_{t=0} \frac{q^{\frac{1}{4}(r+2-t-2s)r_1 + \frac{1}{2}(r+2)t + \frac{1}{2}(r+1)(s-4)}\,\Gamma_q(-1-r_1)}{[s]!\,[t]!\,[r+2-s-t]!\,\Gamma_q(r-r_1+2-s)} \times$$

$$\times D_z^{r+2-t}\,D_y^t\,D_x^t \prod_{u=1}^{r+2-s-t} [N_x - t + 1 - u]\, q^{\frac{(2s-r-2)}{4}N_x + \frac{t}{4}N_y + \frac{(t+r+2)}{4}N_z} \tag{58b}$$

As in the classical case [19] the explicit form of this operator may be checked by the intertwining property:

$$_q\mathcal{D}_3^{r+2}\,\Gamma_3(X)_{r_1,r_2} = \Gamma_3(X)_{r_1',r_2'}\, _q\mathcal{D}_3^{r+2}, \qquad r_1' = -r_2 - 2,\ r_2' = -r_1 - 2, \tag{59}$$

where $X = E_{ii} - E_{i+1,i+1}$, $E_{ij}, i \neq j$, $\Gamma_3(X)_{r_1,r_2}$ is from (14) as in (29), (32).

The states in the subrepresentation are given by $v_{\ell k j}$, with $\ell, k, j \in \mathbb{Z}_+$, $k \leq r+1$.

Let us illustrate the above by the simplest classical example of $r = -1$, $q = 1$. The equation (56a) has the form:

$$\left((r_1 - x\partial_x)\,\partial_z - \partial_y\,\partial_x\right) f(x,y,z) = 0. \tag{60}$$

It is easy to see that all solutions of (60) are given by:

$$f_{j\ell}(x,y,z) = x^j\, y^\ell \sum_{n=0}^{\min(j,\ell)} \frac{j!\,\ell!\,\Gamma(r_1 - j + 1)}{n!\,(j-n)!\,(l-n)!\,\Gamma(r_1 - j + 1 + n)} \left(\frac{z}{xy}\right)^n =$$
$$= x^j\, y^\ell\, {}_2F_1(-j, -\ell\, ;\, r_1 - j + 1\, ;\, \frac{z}{xy})\,, \quad j, \ell \in \mathbb{Z}_+\,. \qquad (61)$$

The same states are valid in the $q$-case since one has from (51):

$$v_{\ell 0 j}|_{r=-1} = \sum_{n=0}^{\ell} (-1)^\ell \binom{\ell}{n}_q q^{\frac{1}{4}\{-\ell j + n(r_1 + \ell)\}} \times$$
$$\times \frac{\Gamma_q(r_1 + 1)}{\Gamma_q(r_1 - j + 1)} \frac{\Gamma_q(1 + r_1 - j + \ell)}{\Gamma_q(1 + r_1 - j + n)} \frac{[j]!}{[j-n]!}\, x^{j-n}\, z^n\, y^{\ell-n} =$$
$$= q^{-\frac{1}{4}\ell j} \frac{\Gamma_q(r_1 + 1)}{\Gamma_q(r_1 - j + 1)} \frac{\Gamma_q(j - r_1)}{\Gamma_q(j - r_1 - \ell)} \times$$
$$\times x^j\, y^\ell\, {}_2F_1^q(-j, -\ell\, ;\, r_1 - j + 1\, ;\, q^{\frac{1}{4}(r_1 + \ell)} \frac{z}{xy})\,. \qquad (62)$$

Alternatively one may check that this is the general solution of (58) for $r = -1$.

**3.6.** Finally, we discuss the case when $r_1, r_2 \in \mathbb{Z}$, $r + 1 = r_1 + r_2 + 1 \in \mathbb{Z}_+$, $r_1 r_2 \in -\mathbb{Z}_+$.

We first take the case $r_1 \in \mathbb{Z}_+$, thus $r_2 \in -\mathbb{N}$ and $r < r_1$. Then we restrict $r_2 < -1$. In this situation there is an invariant subspace which is singled out by the *two* equations (28) and (58), however, in (58b) the two factorials $\Gamma_q(-1 - r_1)/\Gamma_q(r - r_1 + 2 - s)$ are replaced by the regular expression: $(-1)^{r+1-s}[r_1 - r - 1 + s]!/[r_1 + 1]!$. The states in the subrepresentation are given by $v_{\ell k j}$, with $\ell, k, j \in \mathbb{Z}_+$, and *either* $k, \ell \geq 0$, $-r_2 - 1 \leq j \leq r_1$, *or* $\ell \geq 0$, $0 \leq k \leq r + 1$, $0 \leq j \leq -r_2 - 2$. (The last statement is a conjecture proved in partial cases.)

Then we take the case $r_1 \in \mathbb{Z}_+$, $r_2 = -1$. In this situation already the classical operator (56) has to be modified. It has the form (up to overall multiplicative non-zero constant) [19]:

$$\mathcal{D}_3^{r+2} = (\hat{X}_2^+)^{r+2}\,(\hat{X}_1^+)^{r+2}\,, \qquad (63)$$

which can be obtained from (56) by multiplying (56b) with $r_2 + 1$ and then taking the limit $r_2 \to -1$. Thus, it is a composition of the two operators (28) and (31). Since these operators are invariant, the invariant subspace is given by the solutions of (28) only. The latter conclusion may be reached also from the previous subcase taking $r_2 \to -1$; then the second set becomes empty and we are left with what is claimed here.

We then take the case $r_2 \in \mathbb{Z}_+$, $r_1 < -1$. In this situation there is an invariant subspace which is singled out by the *two* equations (31) and (58), The states in the subrepresentation are given by $v_{\ell k j}$, with $\ell, k, j \in \mathbb{Z}_+$, and *either* $k, j \geq 0$, $-r_1 - 1 \leq \ell \leq r_2$, *or* $j \geq 0$, $0 \leq k \leq r+1$, $0 \leq \ell \leq -r_1 - 2$. (The last statement is a conjecture proved in partial cases.)

Finally we take the case $r_2 \in \mathbb{Z}_+$, $r_1 = -1$. In this situation again the classical operator (56) has to be modified. It has the form [19]:

$$\mathcal{D}_3^{r+2} = (\hat{X}_1^+)^{r+2} (\hat{X}_2^+)^{r+2}, \tag{64}$$

which can be obtained from (56) by multiplying (56b) with $r_1 + 1$ and then taking the limit $r_1 \to -1$. Thus, it is a composition of the two operators (28) and (31). Since these operators are invariant, the invariant subspace is given by the solutions of (31) only. The latter conclusion may be reached also from the previous subcase taking $r_1 \to -1$.

Clearly, if $r_1, r_2, r+1 \notin \mathbb{Z}_+$ the representation in (25) is infinite dimensional and irreducible since by (25a) the powers of $x$ are arbitrary, by $(25f'')$ the powers of $y$ are arbitrary, by $(25f')$ the powers of $z, x, y$ are arbitrary.

## 4. Newton diagrams.

It this Section we give a visualization of the representation spaces. Each state is represented by a point on an integer lattice in $n(n-1)/2$ dimensions, i.e., on $\mathbb{Z}_+^{n(n-1)/2}$. For a finite-dimensional subrepresentation the number of these points is finite and the hull of these points is a convex polyhedron in $\mathbb{R}_+^{n(n-1)/2}$. Such a polyhedron (not necessarily convex) is called a Newton diagram [18]. In the present context this notion was introduced in [17], where also some examples in the case of functions in one and two variables were given (for $q = 1$), when the figures are planar (polygons). Below, we give explicitly the Newton diagrams for $n = 3$. Moreover, we introduce also infinite Newton diagrams to depict the infinite-dimensional nontrivial subrepresentations.

**4.1. Finite Newton diagrams for n=3.** Fix $r_k \in \mathbb{Z}_+$. Then the Newton diagram is given by the points with integer coordinates $j, \ell, k$ in $\mathbb{Z}_+^3$ such that:

$$0 \leq j + k + \ell \leq r, \tag{65a}$$

$$0 \leq j \leq r_1, \tag{65b}$$

$$0 \leq \ell \leq r_2, \tag{65c}$$

cf. the **Figure**. The polyhedron formed by these points is planar only for $r_1 = 0$ or $r_2 = 0$ in which case it is a triangle (only (65a) is relevant since $r = r_2$ or $r = r_1$). (The case $r_1 = 0$ was given in [17].)

Fix a point $j, \ell, k$. This is represented by the state $v_{\ell k j}$. Then, the number of states is:

$$\sum_{j=0}^{r_1} \sum_{k=0}^{r-j} \sum_{\ell=0}^{\min(r-k-j,r_2)} 1 = \sum_{j=0}^{r_1} \sum_{k=0}^{r_1-j} \sum_{\ell=0}^{r_2} 1 + \sum_{j=0}^{r_1} \sum_{k=r_1-j+1}^{r-j} \sum_{\ell=0}^{r-k-j} 1 = \tag{66}$$

$$= \frac{(r_1+1)(r_1+2)(r_2+1)}{2} + \frac{(r_1+1)r_2(r_2+1)}{2} = d_{r_1,r_2} ,$$

as expected (cf. (33)).

Note that such diagrams have an advantage over the usual weight diagrams for $sl(3)$ and $su(3)$ which are degenerate. For instance, consider the adjoint representation obtained for $r_1 = r_2 = 1$. The weight diagram consists of two orbits of the Weyl group, one with six points with multiplicity one, and the other with one point with multiplicity two. To the latter point in our diagram correspond the two states:

$$v_{101} = q^{-\frac{1}{4}} ([2]_q xy - q^{-1} z) , \tag{67a}$$
$$v_{010} = q^{-\frac{1}{4}} ([2]_q z - qxy) , \tag{67b}$$

which are linearly independent.

**4.2. Infinite Newton diagrams for n=3.**  Here either $r_1 \notin \mathbb{Z}_+$ or $r_2 \notin \mathbb{Z}_+$ and the considerations run in parallel with Subsections 3.2, 3.4, 3.5. Below $j, \ell, k \in \mathbb{Z}_+$.

**4.2.1.**  For $r_1 \in \mathbb{Z}_+$ and $r_2, r+1 \notin \mathbb{Z}_+$ the Newton diagram is given by the points with coordinates:

$$0 \leq k , \tag{68a}$$
$$0 \leq j \leq r_1 , \tag{68b}$$
$$0 \leq \ell . \tag{68c}$$

**4.2.2.**  For $r_2 \in \mathbb{Z}_+$ and $r_1, r+1 \notin \mathbb{Z}_+$ the Newton diagram is given by the points with coordinates:

$$0 \leq k , \tag{69a}$$
$$0 \leq j , \tag{69b}$$
$$0 \leq \ell \leq r_2 . \tag{69c}$$

**4.2.3.** For $r+1 \in \mathbb{Z}_+$ and $r_1, r_2 \notin \mathbb{Z}_+$ the Newton diagram is given by the points with coordinates:

$$0 \leq k \leq r+1, \qquad (70a)$$
$$0 \leq j, \qquad (70b)$$
$$0 \leq \ell. \qquad (70c)$$

**4.2.4.** For $r_1, r+1 \in \mathbb{Z}_+$ and $r_2 + 1 \in -\mathbb{N}$ the Newton diagram is given by *two* sets of points with coordinates:

$$0 \leq k, \qquad (71a)$$
$$-r_2 - 1 \leq j \leq r_1, \qquad (71b)$$
$$0 \leq \ell. \qquad (71c)$$
$$0 \leq k \leq r+1, \qquad (72a)$$
$$0 \leq j \leq -r_2 - 2, \qquad (72b)$$
$$0 \leq \ell. \qquad (72c)$$

**4.2.5.** For $r_1 = r+1 \in \mathbb{Z}_+$ and $r_2 = -1$ the Newton diagram is given by (68). It can be obtained formally from the previous case by setting $r_2 = -1$, then (71) coincides with (68), while (72) is empty.

**4.2.6.** For $r_2, r+1 \in \mathbb{Z}_+$ and $r_1 + 1 \in -\mathbb{N}$ the Newton diagram is given by *two* sets of points with coordinates:

$$0 \leq k, \qquad (73a)$$
$$0 \leq j \qquad (73b)$$
$$-r_1 - 1 \leq j \leq r_2,. \qquad (73c)$$
$$0 \leq k \leq r+1, \qquad (74a)$$
$$0 \leq \ell, \qquad (74b)$$
$$0 \leq j \leq -r_2 - 2. \qquad (74c)$$

**4.2.7.** For $r_2 = r+1 \in \mathbb{Z}_+$ and $r_1 = -1$ the Newton diagram is given by (69). It can be obtained formally from the previous case by setting $r_1 = -1$, then (73) coincides with (69), while (74) is empty.

**Acknowledgments.**

L.C.B. was supported in part by the U. S. Dept. of Energy under contract # DE-SG03-93ER40757. V.K.D. would like to thank the INFN, Sezione di Genova, for supporting two visits to Genova during which part of this work was done.

**Figure Caption.**

Newton diagram for the finite-dimensional representations of $U_q(sl(3))$.

**Figure**

Newton diagram for the finite-dimensional representations of $\mathbf{U}_q(\mathrm{sl}(\mathbf{3}))$.